\documentstyle[psfig,amssymb]{llncs}

\begin{document}

\newcommand{\eps}{\epsilon}
\newcommand{\f}{{\rm 1st}}
\newcommand{\s}{{\rm 2nd}}
\newcommand{\bof}{{\bf 1st}}
\newcommand{\bos}{{\bf 2nd}}
\newcommand{\boH}{{\bf H}}
\newcommand{\boV}{{\bf V}}
\newcommand{\Neg}{\!-\!}

\title{Who Wins Domineering on Rectangular Boards?}

\author{Michael Lachmann \inst{1}, Cristopher Moore \inst{1,2}, and
Ivan Rapaport \inst{3}}

\institute{
Santa Fe Institute, 1399 Hyde Park Road, Santa Fe NM 87501 
{\tt dirk@santafe.edu}
\and
Computer Science Department and Department of Physics and Astronomy, 
University of New Mexico,
Albuquerque NM 87131 
{\tt moore@cs.unm.edu}
\and
Departamento de Ingenieria Matematica, Facultad de Ciencias Fisicas y
Matematicas, Universidad de Chile, Casilla 170/3, Correo 3, Santiago,
Chile {\tt irapapor@dim.uchile.cl}}

\maketitle

\begin{abstract}
Using mostly elementary considerations, we find out who wins the game
of Domineering on all rectangular boards of width 2, 3, 5, and 7.  We
obtain bounds on other boards as well, and prove the existence of
polynomial-time strategies for playing on all boards of width 2, 3, 4,
5, 7, 9, and 11.  We also comment briefly on toroidal and cylindrical
boards.
\end{abstract}

\section{Introduction}

Domineering or Crosscram is a game invented by G\"oran Andersson and
introduced to the public in \cite{gardner}.  Two players, say Vera and
Hepzibah, have vertical and horizontal dominoes respectively.  They
start with a board consisting of some subset of the square lattice and
take turns placing dominoes until one of them can no longer move.  For
instance, the $2 \times 2$ board is a win for the first player, since
whoever places a domino there makes another space for herself while
blocking the other player's moves.

A beautiful theory of combinatorial games of this kind, where both
players have perfect information, is expounded in \cite{conway,bcg}.
Much of its power comes from dividing a game into smaller subgames,
where a player has to choose which subgame to make a move in.  Such a
combination is called a {\em disjunctive sum}.  In Domineering this
happens by dividing the remaining space into several components, so
that each player must choose in which component to place a domino.

Each game is either a win for Vera, regardless of who goes first, or
Hepzibah regardless of who goes first, or the first player regardless
of who it is, or the second regardless of who it is.  These correspond
to a value $G$ which is positive, negative, fuzzy, or zero, i.e.\ $G >
0$, $G < 0$, $G \,\|\, 0$, or $G = 0$.  (By convention Vera and
Hepzibah are the left and right players, and wins for them are
positive and negative respectively.)  However, we will often
abbreviate these values as $G=V$, $H$, $\f$, or $\s$.  We hope this
will not confuse the reader too much.

In this paper, we find who wins Domineering on all rectangles,
cylinders, and tori of width 2, 3, 5, and 7.  We also obtain bounds on
boards of width 4, 7, 9, and 11, and partial results on many others.
We also comment briefly on toroidal and cylindrical boards.

Note that this is a much coarser question than calculating the actual
game-theoretic values of these boards, which determine how they act
when disjunctively summed with other games.  Berlekamp \cite{berl} has
found exact values for $2 \times n$ rectangles with $n$ odd, and
approximate values to within an infinitesimal or `ish' (which
unfortunately can change who wins in unusual situations) for other
positions of width 2 and 3.  In terms of who wins, the $8 \times 8$
board and many other small boards were recently solved by Breuker,
Uiterwijk and van den Herik using a computer search with a good system
of transposition tables in \cite{buh,uiter,recent}.  We make use of
these results below.

\section{$2 \times n$ boards}

On boards of width 2, it is natural to consider dividing it into two
smaller boards of width 2.  At first glance, Vera (the vertical
player) has greater power, since she can choose where to do this.
However, she can only take full advantage of this if she goes first.
Hepzibah (the horizontal player) has a greater power, since whether
she goes first or second, she can divide a game into two simply by
{\em not} placing a domino across their boundary.  We will see that,
for sufficiently large $n$, this gives Hepzibah the upper hand.

We will abbreivate the value of the $2 \times n$ game as $[n]$.

Let's look at what happens when Vera goes first, and divides a board
of length $m+n+1$ into one of length $m$ and one of length $n$.  If
she can win on both these games, i.e.\ if $[m]=[n]=V$, clearly she
wins.  If $[m]=V$ and $[n]=\s$, Hepzibah will
eventually lose in $[m]$ and be forced to play in $[n]$, whereupon Vera replies there and wins.  Finally, if $[m]=[n]=\s$, Vera replies to Hepzibah in both and wins.  Since
Vera can win if she goes first, $[m+n+1]$ must be a win
either for the first player or for V.  This gives us the following
table for combining boards of lengths $m$ and $n$ into boards of
length $m+n+1$:
\begin{equation} \begin{array}{r|cc}
[m+n+1] & \s & V \\ \hline
     \s & \f \mbox{ or } V & \f \mbox{ or } V \\
      V & \f \mbox{ or } V & \f \mbox{ or } V 
\end{array} \label{tab1} \end{equation}
This table can be summarized by the equation
\begin{equation}
\mbox{If $[m] \ge 0$ and $[n] \ge 0$, then $[m+n+1] \,\,\|\!\!> 0$.}
\label{eqn1}
\end{equation}

Hepzibah has a similar set of tools at her disposal.  By declining to
ever place a domino across their boundaries, she can effectively play
$[m+n]$ as a sum of $[m]$ and $[n]$ for whichever $m$ and $n$ are the
most convenient.  If Hepzibah goes first, she can win whenever
$[m]=\f$ and either $[n]=\s$ or $[n]=H$, by playing first in $[m]$ and
replying to Vera in $[n]$.  If $[m]=[n]=H$ she wins whether she goes
first or second, and if $[m]=\s$ and $[n]=H$, the same is true since
she plays in $[n]$ and replies to Vera in $[m]$.  Finally, if
$[m]=[n]=\s$, she can win if she goes second by replying to Vera in
both games.  This gives the table
\begin{equation} \begin{array}{r|ccc}
[m+n] & \f & \s & H \\ \hline 
   \f & ?                & \f \mbox{ or } H & \f \mbox{ or } H \\
   \s & \f \mbox{ or } H & \s \mbox{ or } H & H \\
    H & \f \mbox{ or } H & H & H
\end{array} \label{tab2} \end{equation}
which can be summarized by the equation
\begin{equation} 
[m+n] \le [m] + [n] \label{eqn2}.
\end{equation}
This simply states that refusing to play across a vertical boundary
can only make things harder for Hepzibah.

These two tables alone, in conjunction with some search by hand and by
computer, allow us to determine the following values.  Values derived
from smaller games using Tables \ref{tab1} and \ref{tab2} are shown in
plain, while those found in other ways, such as David Wolfe's
Gamesman's Toolkit \cite{wolfe}, our own search program, or
Berlekamp's solution for odd lengths \cite{berl} are shown in bold.
\[ \begin{array}{lcclcclcclc}
0 & \bos & \hspace{5mm} 
           & 10 & \f   & \hspace{5mm} 
                         & 20 & H    & \hspace{5mm} 
                                       & 30 & H \\
1 & \boV & & 11 & \f   & & 21 & H    & & 31 & \boH \\
2 & \bof & & 12 & H    & & 22 & \boH & & 32 & H \\
3 & \bof & & 13 & \bos & & 23 & \f   & & 33 & H \\
4 & \boH & & 14 & \f   & & 24 & H    & & 34 & H \\
5 & \boV & & 15 & \f   & & 25 & H    & & 35 & H \\
6 & \f   & & 16 & H    & & 26 & H    & & 36 & H \\
7 & \f   & & 17 & H    & & 27 & \f   & & 37 & H \\
8 & H    & & 18 & \bof & & 28 & H    & & 38 & H \\
9 & \boV & & 19 & \f   & & 29 & H    & & 39 & H
\end{array} \] 
In fact, $[n]$ is a win for Hepzibah for all $n \ge 28$.

Some discussion is in order.  Once we know that $[4]=H$, we have
$[4k]=H$ for all $k \ge 1$ by Table \ref{tab2}.  Combining Tables
\ref{tab1} and \ref{tab2} gives $[6]=[7]=\f$, since these are both
$\f$ or $V$ and $\f$ or $H$.  A similar argument gives $[10]=[11]=\f$
and $[14]=[15]=\f$, once we learn through search that $[13]=\s$ (which
is rather surprising, and breaks an apparent periodicity of order 4).

Combining $[13]$ with multiples of 4 and with itself gives
$[13+4k]=[26+4k]=[39+4k]=H$ for $k \ge 1$.  Since $26=24+2=13+13$, we
have $[26]=H$ since it is both $\f$ or $H$ and $\s$ or $H$, giving
$[39]=H$ since 39=26+13.

Since $19=9+9+1=17+2$, we have $[19]=\f$ since it is both $\f$ or $V$
and $\f$ or $H$.  Similarly $23=9+13+1=21+2$ and $27=13+13+1=25+2$ so
$[23]=[27]=\f$.  A computer search gives $[22]=H$, and since
$35=22+13$ we have $[35]=H$.

So far, we have gotten away without using the real power of game
theory.  However, for $[31]$ we have found no elementary proof, and it
is too large for our search program.  Therefore, we turn to
Berlekamp's beautiful solution for $2 \times n$ Domineering when $n$
is odd (Ref. \cite{berl}), evaluate it with the Gamesman's Toolkit
\cite{wolfe}, and find the following (see Refs. \cite{conway,bcg,berl}
for notation):
\begin{equation} 
[31] \,=\, \frac12 - 15 \cdot \left( \frac14 + \int^{3/4} * \right) 
           + \int^{3/4} \int_{1/2}^{1/2*} 3\frac78 \,=\, 
\left\{ 2 \,|\, 0 \,\|\, {-\frac12} \,|\, {-2} \right\} |\, {-\frac52} 
 \,<\, 0
\end{equation}
Thus $[31]$ is negative and a win for Hepzibah.  This closes the last
loophole, telling us who wins the $2 \times n$ game for all $n$.

\section{Boards of width 3, 4, 5, 7, 9, 11 and others}

The situation for rectangles of width 3 is much simpler.  While
Equation~\ref{eqn1} no longer holds since Vera cannot divide the board
in two with her first move, Equation~\ref{eqn2} holds for all widths,
since Hepzibah can choose not to cross a vertical boundary between two
games.  A quick search shows that $[3 \times n]=H$ for $n=4,5,6$ and
$7$, so we have for width 3
\[ \mbox{$[1]=V$, $[2]=[3]=\f$, and $[n]=H$ for all $n \ge 4$.} \]
For width 5, we obtain
\[ \mbox{$[1]=[3]=V$, $[2]=[4]=H$, $[5]=\s$, and $[n]=H$ for
all $n \ge 6$.} \]
For width 7, Breuker, Uiterwijk and van den Herik found by computer
search (Refs. \cite{buh,uiter,recent}) that $[4]=[6]=[9]=[11]=H$.
Then $[8]=[10]=H$, and combining this with searches on small boards we
have
\[ \mbox{$[1]=[3]=[5]=V$, $[2]=[7]=\f$, $[4]=[6]=H$, and $[n]=H$ 
for all $n \ge 8$.} \]

In all these cases, we were lucky enough that $[n]=\s$ or $H$ for
enough small $n$ to generate all larger $n$ by addition.  This becomes
progressively rarer for larger widths.  However, we have some partial
results on other widths.  For width 4, Uiterwijk and van den Herik
(Refs. \cite{uiter,recent}) found by computer search that
$[8]=[10]=[12]=[14]=H$, so $[n]=H$ for all even $n \ge 8$.  They also
found that $[15]=[17]=H$, so
\[ [4 \times n]=H \mbox{ for all } n \ge 22. \]
This leaves $[4 \times 19]$ and $[4 \times 21]$ as the only unsolved
boards of width 4.

As a general method, whenever we can find a length for which Hepzibah
wins by some positive number of moves (rather than by an
infinitesimal), then she wins on any board long enough to contain a
sufficient number of copies of this one to overcome whatever advantage
Vera might have on smaller boards.  Game-theoretically, if $[n] < -r$,
then $[m] < 0$ whenever $m \ge (1/r) \max_{\,l < n} \,[l]$.

For width 9, for instance, we have $[1]=4$, $[2]=\frac32 \,|\, 0
\,\|\, {-\frac12} \,|\, {-\frac52}$, $[3]=5 \,|\, 3 \,\|\,
\frac{11}{4} \,|\, \frac14$, and $[4] \le [2]+[2] = 1 \,|\, {-\frac12}
\,\|\, {-1} \,|\, {-\frac52} < {-\frac12}$.  By summing these, it is
easy to show that
\[ [9 \times n]=H \mbox{ for all } n \ge 22. \]
Similarly, for width 11 we have $[1]=5$ and $[2]=1 \,|\, \left\{
\frac12 \,|\, {-1} \,\|\, {-\frac32} \,|\, {-\frac72} \right\}$.  Then
$[8] \le [2]+[2]+[2]+[2] = 1 \,|\, {-\frac12} \,\|\, {-1} \,|\,
{-\frac52} < {-\frac12}$ and $[16] \le [8]+[8] \le {-\frac32}$, so
\[ [11 \times n]=H \mbox{ for all } n \ge 56. \]
Unfortunately, for all other widths greater than 7, either [2] or
$[2]+[2]$ is positive and $[3]$ is as well, so without some way to
calculate values for $[4]$ or more we can't establish this kind of
bound.  Nor do we know of any length for which Hepzibah wins on width
8, or a proof that she wins any board of width 6 by a positive amount.

To get results on boards of other widths, we can use a variety of
tricks.  First of all, just as Hepzibah can choose to cross a vertical
boundary between games, Vera can choose not to cross a horizontal one.
Thus Equation~\ref{eqn2} is one of a dual pair,
\begin{eqnarray}
[m \times (n_1+n_2)] & \le & [m \times n_1] + [m \times n_2] \\ 
{}[(m_1+m_2) \times n] & \ge & [m_1 \times n] + [m_2 \times n]
\end{eqnarray}
Another useful rule is that $[n \times n]=\f$ or $\s$, since neither
player can have an advantage on a square board.  In fact, in
game-theoretic terms $[n \times n]+[n \times n]=0$, so if Hepzibah
goes second she can win by mimicking Vera's move, rotated $90^\circ$,
in the other board.  More generally we have
\[ \begin{array}{l}
     \mbox{If $[n \times n]=\f$, then } [n \times kn]=
       \left\{ \begin{array}{ll}
         \s \mbox{ or } H & \mbox{ for even } k > 1 \\
         \f \mbox{ or } H & \mbox{ for odd } k > 1
       \end{array} \right. \\
     \mbox{If $[n \times n]=\s$, then $[n \times kn]=\s$ or $H$ 
       for all $k>1$.}
\end{array} \]
For instance, this tells us that $[6 \times 12]=\s$ or $H$, and since
$[6 \times 4]=\f$ and $[6 \times 8]=H$ (Ref. \cite{buh}) we also have
$[6 \times 12]=\f$ or $H$.  Therefore $[6 \times 12]=H$ and
\[ [6 \times (4+4k)]=H \mbox{ for all } k \ge 1. \]
We can also use our addition rules backward; since no two games can
sum to a square in a way that gives an advantage to either player, 
\begin{equation}
\mbox{If $m < n$ and $[m \times n]=\s$ or $V$, then $[(n-m) \times n]
\ne V$}
\end{equation}
and similarly for the dual version.  

\begin{figure}
\begin{center}
\psfig{figure=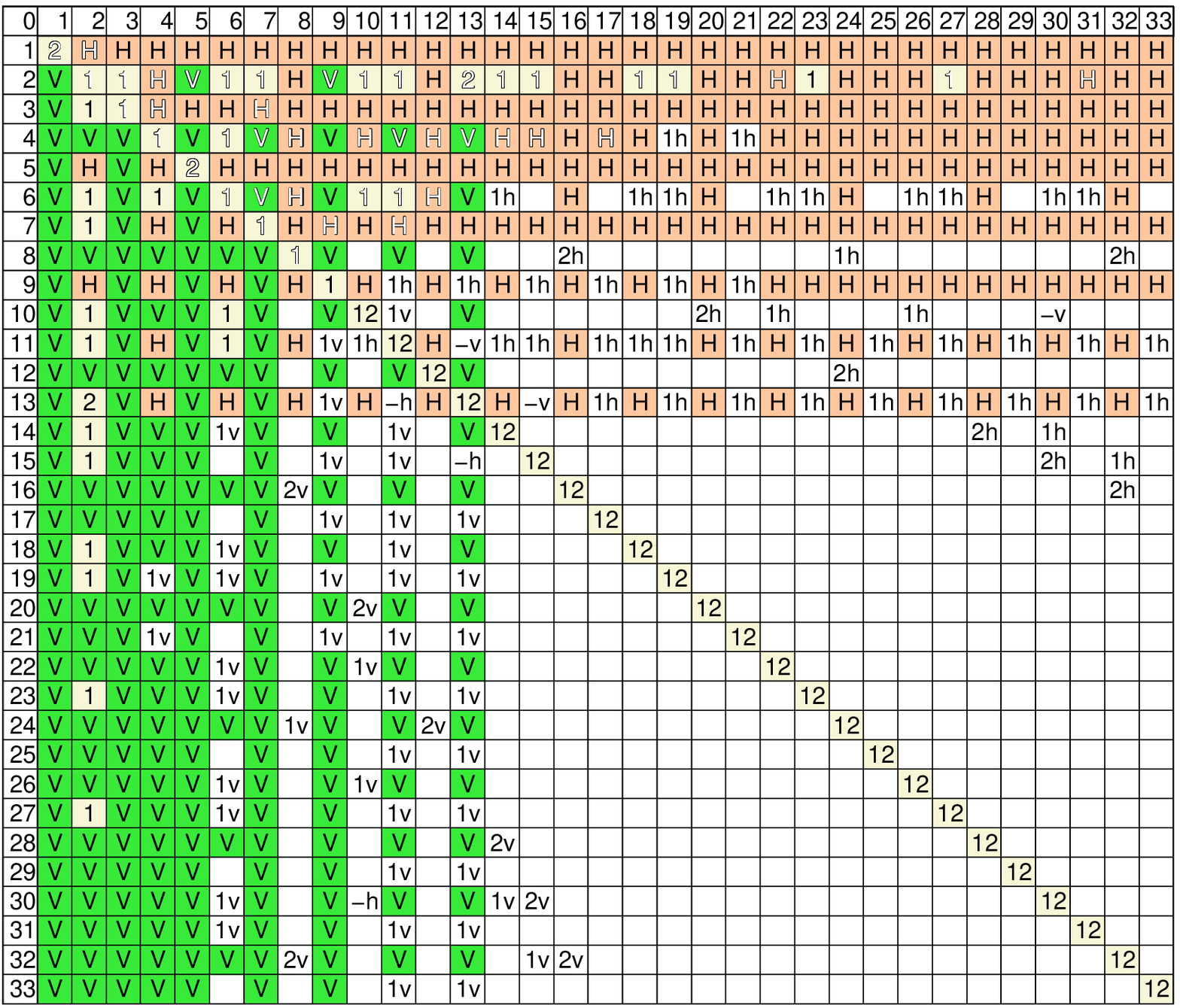,width=5in}
\end{center}
\caption{What we know so far about who wins Domineering on rectangular
boards.  $1$, $2$, $H$ and $V$ mean a win for the first player, second
player, Hepzibah and Vera respectively.  Things like ``1h'' mean either
$\f$ or $H$ (i.e.\ all we know is that Hepzibah wins if she goes first)
and ``-v'' means that it is not a win for Vera all the time.  Values
outlined in black are those provided by search or other methods; all
others are derived from these using our rules or by symmetry.}
\label{bigtab}
\end{figure}

Using the results of Refs. \cite{berl} and \cite{buh}, some computer
searches of our own, and a program that propagates these rules as much
as possible gives the table shown in Figure~\ref{bigtab}.  It would be
very nice to deduce who wins on some large squares; the $9 \times 9$
square is the largest known so far (Ref. \cite{recent}).  We note that
if $[9 \times 13]=\f$ then $[13 \times 13]=\f$ since $[4 \times
13]=V$.

\section{Playing on cylinders and tori}

\begin{figure}
\centerline{\psfig{figure=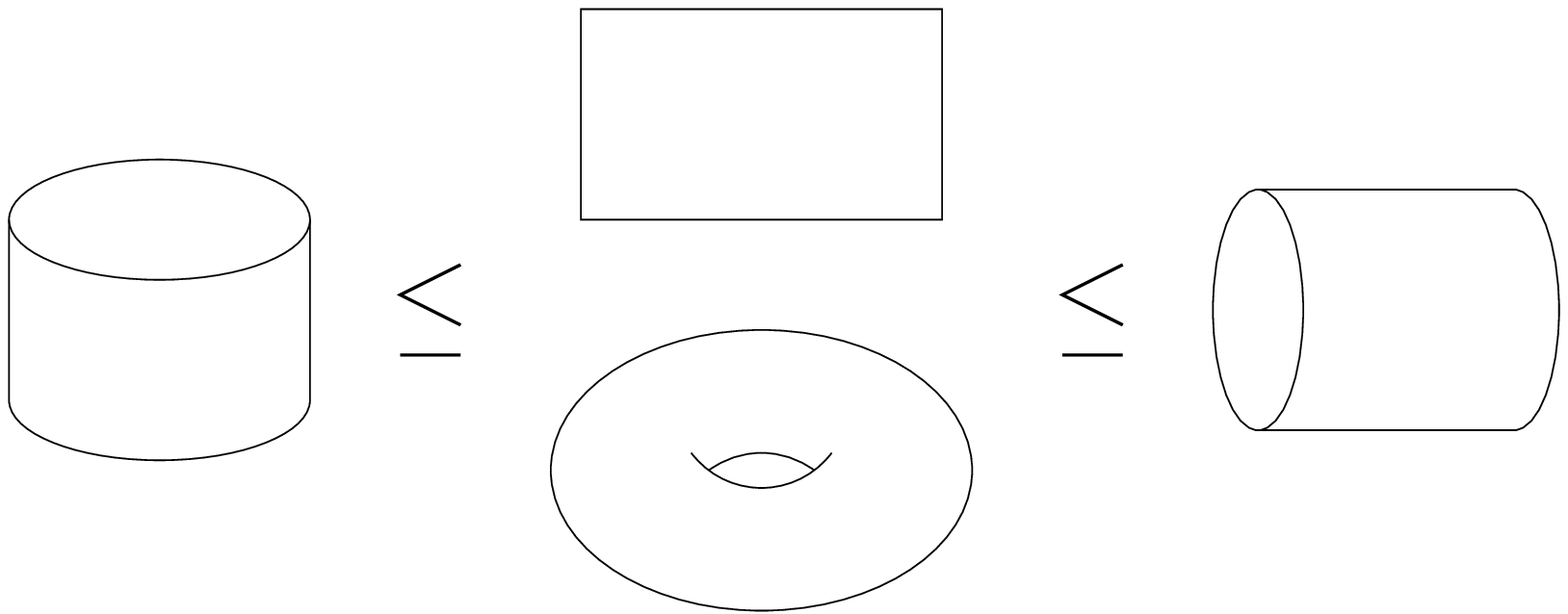,width=3in}}
\caption{Inequalities between rectangular, cylinder, and toroidal
boards of the same size.  Cutting vertically can only hurt Hepzibah,
while cutting horizontally can only hurt Vera.}
\label{torusfig}
\end{figure}

On a torus, Hepzibah can choose not to play across a vertical boundary
and Vera can choose not to play across a horizontal one.  Thus cutting
a torus, or pasting a rectangle along one pair of edges, to make a
horizontal or vertical cylinder gives the inequalities shown in
Figure~\ref{torusfig}.  Note that there is no obvious relation between
the value of a rectangle and that of a torus of the same size.

While it is easy to find who wins on tori and cylinders of various
small widths, we do the analysis here only for tori of width 2.  Since
Vera's move takes both squares in the same column, these boards are
equivalent to horizontal cylinders like those shown on the left of
Fig.~\ref{torusfig} (or, for that matter, M\"obius strips or Klein
bottles).  Therefore, if Hepzibah can win on the rectangle of length
$n$, she can win here as well.

The second player has slightly more power here than she did in on the
rectangle, since the first player has no control over the effect of
her move.  If Vera goes first, she simply converts a torus of length
$n$ into a rectangle of length $n-1$, and if Hepzibah goes first, Vera
can choose where to put the rectangle's vertical boundary, in essence
choosing Hepzibah's first move for her.  On the other hand, in the
latter case Hepzibah gets to play again, and can treat the remainder
of the game as the sum of two rectangles and a horizontal space.

These observations give the following table for tori of width 2:
\[ \begin{array}{cc|c}
[n]_{\rm rect} & [n-1]_{\rm rect} & [n]_{\rm torus} \\ \hline
H & & H \\
\f & \f \mbox{ or } H & H \\
\f & \s \mbox{ or } V & \f \\
\s \mbox{ or } V & \f \mbox{ or } H & \s \mbox{ or } H 
\end{array} \]
These and our table for $2 \times n$ rectangles determine $[n]_{\rm
torus}$ for all $n$ except 5, 9, and 13.  Vera loses all of these if
she plays first, since she reduces the board to a rectangle which is a
win for Hepzibah.  For 5 and 9, Vera wins if Hepzibah plays first by
playing in such a way that Hepzibah's domino is in the center of the
resulting rectangle, creating a position which has zero value.  Thus
these boards are wins for the $\s$ player.  For 13, Hepzibah wins
since (as the Toolkit tells us) all of Vera's replies to Hepzibah
leave us in a negative position.

Our computer searches show that $n \times n$ tori are wins for the
$\s$ player when $n=1$, $3$, or $5$, and for the $\f$ player when
$n=2$, $4$, or $6$.  We conjecture that this alternation continues,
and that square tori of odd and even size are wins for the $\s$ and
$\f$ players respectively.  We note that a similar argument can be
used to show that 9 is prime.

\section{Polynomial-time strategies}

While correctly playing the sum of many games is PSPACE-complete in
general \cite{morris}, the kinds of sums we have considered here are
especially easy to play.  For instance, if $[m]$ and $[n]$ are both
wins for Hepzibah, she can win on $[m+n]$ by playing wherever she
likes if she goes first, and replying to Vera in whichever game Vera
chooses thereafter.  Thus if we have strategies for both these games,
we have a strategy for their sum.  All our additive rules are of this
kind.

Above, we showed for a number of widths that boards of any length can
be reduced to sums of a finite number of lengths.  Since each of these
can be won with some finite strategy, and since sums of them can be
played in a simple way, we have proved the following theorem:

\begin{theorem}
For boards of width 2, 3, 4, 5, 7, 9, and 11, there exist
polynomial-time strategies for playing on boards of any length.
\end{theorem}

Note that we are not asking that the strategy produce optimum play, in
which Hepzibah (or on small boards, $\f$, $\s$ or Vera) wins by as
many moves as possible, but only that it tells her how to win.

In fact, we conjecture that this theorem is true for boards of any
width.  This would follow if for any $m$ there exists an $n$ such that
Hepzibah wins by some positive number of moves, which in turn implies
that there is some $n'$ such that she wins on all boards longer than
$n'$.  A similar conjecture is made in \cite{uiter}.  Note, however,
that this is not the same as saying that there is a single
polynomial-time strategy for playing on boards of any size.  The size
or running time of the strategy could grow exponentially in $m$, even
if it grows polynomially when $m$ is held constant.

{\bf Acknowledgements.}  We thank Elwyn Berlekamp, Aviezri Fraenkel,
and David Wolfe for helpful communications, and Jos Uiterwijk for
sharing his group's recent results.  I.R.\ also thanks the Santa Fe
Institute for hosting his visit, and FONDECYT 1990616 for their
support.  Finally, C.M. thanks Molly Rose and Spootie the Cat.

\end{document}